\documentclass[letterpaper, 10 pt, conference]{cssconf}  

\IEEEoverridecommandlockouts                              
\usepackage{graphicx}
\usepackage{psfrag}
\usepackage{amsmath} 
\usepackage{amssymb}  

\usepackage{algorithm}
\usepackage{algorithmic}

\newcommand{\RETURN}{\textbf{return: }}

\newcommand{\R}{{\mathbb R}}
\newcommand{\C}{{\mathbb C}}

\newcommand{\VV}{{\mathcal V}}
\newcommand{\BB}{{\mathcal B}}

\def\argmin{\mathop{\rm arg\,min}}

\newtheorem{theorem}{Theorem}
\newtheorem{proposition}{Proposition}

\title{\LARGE \bf
An algebraic geometry approach to\\
nonlinear parametric optimization in control}

\author{Ioannis A. Fotiou*, Philipp Rostalski*, Bernd Sturmfels$^\dagger$ and Manfred Morari*
\thanks{$^*$ Automatic Control Laboratory, Swiss Federal Institute of
Technology, CH-8092 Zurich, Switzerland.}%
\thanks{$^\dagger$ Department of Mathematics, University of California at
Berkeley, 94720-3840 CA, USA.}%
\thanks{Corresponding Author: Email: \texttt{fotiou@control.ee.ethz.ch}}%
}

\begin{document}

\maketitle
\thispagestyle{empty}
\pagestyle{empty}

\begin{abstract}
We present a method for nonlinear parametric optimization based on
algebraic geometry. The problem to be studied, which arises in
optimal control, is to minimize a polynomial function with
parameters subject to semialgebraic constraints. The method uses
Gr\"obner bases computation in conjunction with the eigenvalue
method for solving systems of polynomial equations. In this way,
certain companion matrices are constructed off-line. Then, given
the parameter value, an on-line algorithm is used to efficiently
obtain the optimizer of the original optimization problem in real
time.
\end{abstract}


\section{INTRODUCTION}

Optimal control is a very active area of research with broad
industrial applications \cite{QiBa99}. It is among the few control
methodologies providing a systematic way to perform nonlinear
control synthesis that handles also system constraints. To a great
extent, it is thanks to this capability of dealing with
constraints that model predictive control (MPC) has proven to be
very successful in practice \cite{MoLe99}, \cite{GaPM89}.

Model predictive control uses optimization on-line to obtain the
solution of the optimal control problem in real time. This method
has been proven most effective for applications. Typically, the
optimal control problem can be formulated into a discrete time
mathematical program, whose solution yields a sequence of control
moves. Out of these control moves only the first is applied,
according to the receding horizon control (RHC) scheme.

The optimal control problem is formulated as a mathe\-matical
program, which can be a linear program (LP), a quadratic program
(QP) or a general nonlinear program (NLP). For hybrid systems, the
corresponding mathematical programs can be mixed integer programs
- MILPs, MIQPs or MINLPs \cite{BeMo99}. The class of the
optimization problem depends on the objective function and the
class of systems one wants to derive an optimal controller for.

Technology and cost factors, however, make the implementation of
receding horizon control difficult if not, in some cases,
impossible. To circumvent these issues, the solution of the
optimal control problem is computed off-line, by solving the
corresponding mathematical program parametrically \cite{BMDP02}.
That is, we compute the explicit formula giving the solution of
the program (control inputs) as a function of the problem
parameters (measured state). The solution then is efficiently
implemented on-line as a lookup table.

In the present work, we extend the concept of the explicit
solution to the class of nonlinear polynomial systems with
polynomial cost function. By polynomial systems we mean those
systems, whose state update equation is given by a polynomial
vector field. For this class of systems, the resulting
mathematical program is a nonlinear (polynomial) parametric
optimization problem.

While the \textit{explicit solution} is not generally possible in
the nonlinear case, we stress the fact that a \textit{partial
precomputation} of the optimal control law is still feasible using
algebraic techniques \cite{FoPM05}. In this paper, we use the
eigenvalue method \cite{CoLO98} in conjunction with Gr\"obner
bases computation to perform nonlinear parametric optimization of
polynomial functions subject to polynomial constraints.

\section{PARAMETRIC OPTIMIZATION}
\label{Seq:Parametric}

Let $u \in \R^m$ be the decision-variable vector and $x \in \R^n$
be the parameter vector. The class of optimization problems that
this paper deals with can generally assume the following form:
\begin{equation}
\begin{array}{c}
        \min_{u} \limits \, J(u, x) \qquad
        \mbox{s.t.} \quad g(u,x) \leq 0,
    \end{array}
    \label{Eqn:paramopt}
\end{equation}
where $J(u,x) \in \R[x_1,\ldots,x_n,u_1,\ldots,u_m]$ is the
objective function and $g\in$
$\R[x_1,\ldots,x_n,u_1,\ldots,u_m]^q$ is a vector polynomial
function representing the constraints of the problem. By
parametric optimization, we mean minimizing the function $J(u,x)$
with respect to $u$ for any given value of the parameter $x \in
\mathcal{X} \subseteq \R^n$, where $\mathcal{X}$ is the set of admissible
parameters. Therefore, the polynomial parametric optimization
problem is finding a computational procedure for evaluating the
maps
\begin{equation}
\begin{array}{llcl}
u^*(x):     &\R^n    &\longrightarrow     &\R^m \\
            &x       &\longmapsto         &u^*\\
            \\
J^*(x):     &\R^n    &\longrightarrow     &\R \\
            &x       &\longmapsto         &J^*,\\
\end{array}
\label{Eqn:maps}
\end{equation}
where
\begin{equation}
\begin{array}{lcl}
u^* &=& \argmin_u \limits \, J(u,x) \\
J^* &=& \min_{u} \limits \, J(u, x).
\end{array}
\label{stars}
\end{equation}
For the sake of simplicity, we assume that the feasible set
defined by $g(u,x)$ is compact, therefore the minimum is attained.
Also, in order for (\ref{Eqn:maps}) not to be point-to-set maps, we
focus our attention to one (any) optimizer.

\subsection{Posing the problem}
Our point of departure is the observation that the cornerstone of
continuous constrained optimization are the Karush-Kuhn-Tucker
(KKT) conditions. All local and global minima for problem
(\ref{Eqn:paramopt}) (satisfying certain constraint
qualifications) occur at the so-called ``critical points"
\cite{BoVa04}, namely the solution set of the following system:
\begin{equation}
\begin{array}{rll}
\nabla_u{J(u,x)} + \sum_{i=1}^{q}{\mu_i \nabla_u g_i(u,x)} &=&
0 \\
\mu_i g_i(u,x) &=& 0 \\
\mu_i &\geq& 0 \\
g(u,x) &\leq& 0.
\end{array}
 \label{Eqn:KKT}
\end{equation}
For the class of problems we consider, the two first relations of
the KKT conditions~(\ref{Eqn:KKT}) form a \textit{square system of
polynomial equations}. Various methods have been proposed in the
literature for solving systems of polynomial equations, both
numerical and symbolic \cite{Stur02}, \cite{PaSt00},
\cite{Mano94}. Here we consider symbolic methods since our aim is
to solve the optimization problem parametrically. We should point
out that the underlying philosophy is that we aim at moving as
much as possible of the computational burden of solving the
nonlinear program (\ref{Eqn:paramopt}) off-line, leaving an easy
task for the on-line implementation.

\subsection{Off-line vs. on-line computations}
The explicit representation of the optimal control law as a state
feedback has been successfully investigated for the linear,
quadratic and piecewise affine case. Among other advantages of the
explicit representation is that one is able to analyze the
controller, derive Lyapunov functions \cite{ChrEtal:cdc:04},
perform dynamic programming iterations
\cite{BorEtal:2005:IFA_2217} in an effective way, even compute the
infinite horizon solution for certain classes of constrained
optimal control problems \cite{BaoEtal:cdc:03}.

Unfortunately, such an explicit representation is not always
possible. The enabling factor in the case of linear systems (or
piecewise affine systems) is the fact that the KKT system
(\ref{Eqn:KKT}) can be solved analytically. In the general
polynomial case studied here, we have to solve a system of
(nonlinear) polynomial equations. The next best alternative then
to an explicit solution is to bring the system in such a form, so
that once the parameters are specified, the solution can be
extracted easily and fast.

\section{THE EIGENVALUE METHOD}
\label{Sec:algebraic_techniques}

In this section we briefly describe the method of eigenvalues
(\cite{CoLO98}, Chapter 2, \S 4) for solving systems of polynomial
equations. This method is used in conjunction with Gr\"obner bases
to perform parametric optimization.

\subsection{Solving systems of polynomial equations}
Suppose we have a system of $m$ polynomial equations $f_i$ in $m$
variables $u_i$
\begin{equation}
\begin{array}{lcl}
f_1(u_1,\ldots,u_m) &=& 0 \\
&\cdots& \\
f_m(u_1,\ldots,u_m) &=& 0.
\end{array}
\label{eqsys}
\end{equation}
These equations form an ideal $I \in K[u_1,\ldots,u_m]$, where $K$
denotes an arbitrary field:
\begin{equation}
I := \langle f_1, \ldots, f_m \rangle \; . \label{Eqn:ideal}
\end{equation}
The solution points we are interested in are the points on the
variety over the algebraic closure $\overline{K}$ of $K$,
\begin{equation}
V(I)=\lbrace s \in \overline{K}^m :\; f_1(s)=0, \ldots, f_m(s)=0
\rbrace, \label{Eqn:Variety}
\end{equation}
i.e. the set of common zeros of all polynomials in the ideal $I$.
These points can be computed by means of Gr\"obner bases. An
obvious choice would be a projection-based algorithm by means of
lexicographic Gr\"obner bases, see~(\cite{CoLO92}, Chapter 2, \S
8). Since the computation of a lexicographic Gr\"obner basis is
very time consuming, we focus on a different method.

The first step we take towards solving (\ref{eqsys}) is computing a
Gr\"obner basis with an arbitrary term-order, e.g. graded reverse
lexicographic term-order. We define $G
=\{\gamma_1,\ldots,\gamma_t\}$ to be this Gr\"obner basis of $I$.
\subsection{The generalized companion matrix}
Consider a polynomial function $h \in K[u_1,\ldots,u_m]$. The
Gr\"obner basis $G$ and the division algorithm make it possible to
uniquely write any polynomial $h \in K[u_1,\ldots,u_m]$ in the
following form:
\begin{equation}
h = c_1(u) \gamma_1 + \cdots + c_t(u) \gamma_t + \overline{h}^{G},
\end{equation}
where $\bar{h}^{G}$ is the unique remainder of the division of $h$
with respect to the Gr\"obner basis $G$. The polynomial $h$ can in
turn be multiplied with another polynomial function $f \in
K[u_1,\ldots,u_m]$ and their product expressed as follows:
\begin{equation}
f \cdot h = d_1(u) \gamma_1 + \cdots + d_t(u) \gamma_t + \overline{f \cdot
h}^{G}.
\end{equation}
In the generic case, the ideal $I$ will be
\textit{zero-dimensional}, which means that the corresponding
\textit{quotient ring}
\begin{align}
A = K[u_1,\ldots,u_m] / I \label{qring}
\end{align}
is a finite-dimensional $K$-vector space (\cite{CoLO92}, Chapter
5, \S 2). The quotient ring of an ideal can be thought of as the
set of all polynomials that do not belong to the ideal but belong
to the underlying ring. Denote with $b = [b_1,\ldots,b_l]^T$ the
vector of the \textit{standard monomials}. A monomial is standard
if it is not divisible by any leading monomial of a polynomial in
the Gr\"obner basis. These standard monomials of $G$ form a basis
\begin{equation}
\mathcal{B} = \{ b_1,\ldots, b_l \} \label{eqn::basis}
\end{equation}
for the $K$-vector space $A$. As a result, every remainder can be
expressed with respect to this basis as an inner product
\begin{align}
r_i = a_i^T \cdot b \; ,
\end{align}
where $a_i \in K^l$. We can now define the map $m_h: \; A
\rightarrow A$ as follows: if $\overline{p}^G \in A$, then
\begin{align}
m_h(\overline{p}^G) := \overline{h\cdot p}^G  =
\overline{\overline{h}^G \cdot \overline{p}^G}^G \in A.
\end{align}
The following proposition holds.
\begin{proposition}
Let $h \in K[u_1,\ldots,u_m]$. Then the map $m_h:A \rightarrow A$
is $K$-linear.
\label{sprop}
\end{proposition}
The proof of Proposition \ref{sprop} can be found in
(\cite{CoLO92}, p. 51). Since $A$ is a finite-dimensional vector
space and the map $m_h$ is linear, its representation with respect
to a basis of this vector space is given by a square matrix $M_h$.
The $l \times l$-matrix $M_h$ is called the generalized companion
matrix.

\subsection{Computing the companion matrix}
To compute the matrix $M_h$, assume that we have the basis $\BB =
\{b_1, \ldots, b_l\}$ consisting of the standard monomials $b_i$
of the Gr\"obner basis $G$. Then, for each one of them, compute
the remainder $r_i$ of the polynomial $h \cdot b_i$ with respect to
the Gr\"obner basis $G$:
\begin{align}
\overline{h \cdot b_i}^G = r_i, \quad \forall \; b_i \in \BB.
\end{align}
All $r_i \in A$ can in turn be expressed as an inner product
\begin{align}
r_i = a_i^T \cdot b
\end{align}
with respect to the basis $\BB$. By collecting all vectors $a_i$
for all basis elements \cite{CoLO98}, we can construct a
representation of the map $m_h$ with respect to basis $\BB$, i.e.
calculate the matrix $M_h$ as follows:
\begin{align}
M_h \equiv [a_{ij}] = \left[
\begin{array}{c}
a_1^T \\
\cdots \\
a_l^T
\end{array}
\right].
\end{align}
Computing the companion matrix is a standard algebraic procedure
implemented in various packages, e.g. in Maple~10.

\subsection{Evaluating polynomial functions on a variety}
Consider a polynomial function $h \in \R[u_1,\ldots,u_m]$. The
amazing fact about the matrix $M_h$ is that the set of its
eigenvalues is exactly the value of $h$ over the variety
$\mathcal{V}(I)$ defined by the ideal $I$. More precisely,
$\mathcal{V}(I)$ is the set of all solution points in complex
$m$-space $\C^m$ of the system (\ref{eqsys}). The following
theorem holds.
\begin{theorem}
Let $I \subset \C[u_1,\ldots,u_m]$ be a zero-dimensional ideal,
let $h \in \C[u_1,\ldots,u_m]$. Then, for $\lambda \in \C$, the
following are equivalent:
\begin{enumerate}
\item $\lambda$ is an eigenvalue of the matrix $M_h$

\item $\lambda$ is a value of the function $h$ on the variety
$\mathcal{V}(I)$.
\end{enumerate}
\end{theorem}
The proof can be found in (\cite{CoLO98}, p. 54).

To obtain the coordinates of the solution set of (\ref{eqsys}), we
evaluate the functions
\begin{equation}
\begin{array}{llcl}
h_1:&        u   &\longmapsto     &u_1 \\
                 &\cdots&              \\
h_m:&        u   &\longmapsto     &u_m
\end{array}
\end{equation}
on the variety $\mathcal{V}(I)$ defined by the ideal $I$, where
$u$ above denotes the vector $(u_1,\ldots,u_m)$. This can be done
by means of the associated companion matrices of the functions
$h_i$. The following theorem taken from (\cite{Stur02} p. 22) is
the basis for the calculation of these point coordinates.
\begin{theorem}
\emph{The complex zeros of the ideal I are the vectors of joint
eigenvalues of the companion matrices} $M_{u_1} \dots M_{u_{m}}$,
\emph{that is,}
\begin{align*}
\VV(I) = &\big \lbrace  (u_1, \dots ,u_{m}) \in \R^{m} :\\
 &\exists v \in \R^{m} \; \forall\, i \; :\; M_{u_i}  v = u_i v
\big\rbrace
\end{align*}
\label{th::eigenvalue}
\end{theorem}
It has to be noted that any vector-valued polynomial function $h:
\R^m \longrightarrow \R$ can be evaluated over a zero-dimensional
variety in the same way.

\section{THE ALGORITHM}
\label{sec::algorithm}

In this section, we present the proposed algorithm, which consists
of two parts: the off-line part, where the generalized companion
matrices for the optimization problem are constructed, and the
on-line part where this precomputed information is used and given
the value of the parameter $x$, the optimal solution is
efficiently extracted.

\subsection{Idea}
Under certain regularity conditions, if $J^*$ (defined in
(\ref{stars})) exists and occurs at an optimizer $u^*$, the KKT
system (\ref{Eqn:KKT}) holds at $u^*$. Consequently, $J^*$ is the
minimum of $J(u,x)$ over the semialgebraic set defined by the KKT
equations and inequalities~(\ref{Eqn:KKT}). These conditions can
be separated in a set of inequalities and a square system of
polynomial equations. The method of eigenvalues for solving
systems of polynomial equations as described in
section~\ref{Sec:algebraic_techniques} can be used for the latter.
This method assumes that the ideal generated by the KKT system
(\ref{Eqn:KKT}) is zero-dimensional.

By ignoring the inequalities, a superset of all critical points is
computed and in a second step, all infeasible points are removed.
Finally, among the feasible candidate points those with the
smallest cost function value have to be found via discrete
optimization. By discrete optimization we mean choosing among a
finite set that point, which yields the smallest objective
function value.

\subsection{Off-line Part}
In $K[u_1, . . . , u_m, \mu_1, . . . , \mu_q]$, where $K$ is the
field of rational functions $\R(x_1\ldots,x_n)$ in the parameter
$x$, we define the KKT ideal
\begin{align}
\!\!\!\!\!I_{KKT} = \langle \nabla_{u}{J(u,x)} +
\sum_{i=1}^{q}{\mu_i \nabla_{u} g_i(u,x)}, \; \mu_i g_i(u,x)
\rangle \! \label{Eqn:gradideal}
\end{align}
containing all the equations within the
KKT-system~(\ref{Eqn:KKT}). All critical points for the
optimization problem~(\ref{Eqn:KKT}) and fixed $x$ are the subset
of real points on the KKT-variety
\begin{align}
\VV_{KKT}^{\R} \subseteq \VV_{KKT}=\VV(I_{KKT}) \; .
\label{eq::fullKKT}
\end{align}
Using the method described in
section~\ref{Sec:algebraic_techniques} we can compute these by
means of the generalized companion matrices.

The algebraic part of the algorithm, i.e the computation of the
companion matrices can be done parametrically. For one thing, one
could use Gr\"obner bases computation for the ideal $I_{KKT}$ and
try to compute the corresponding companion matrices $M_{u_i}$ and
$M_{\mu_i}$ directly. Owing to the structure of the polynomial
equations of the $KKT$-system~(\ref{Eqn:gradideal}), this problem
is very poorly conditioned. The difficulties stem from the fact
that the ideal $I_{KKT}$ is by construction decomposable. It
contains terms like $\mu_i g_i(u,x) $ which lead to a reducible
variety $V(I_{KKT})$.

To overcome this obstacle, we factorize the generators of the
Gr\"obner basis (i.e. the polynomials appearing in relation
(\ref{Eqn:gradideal})) and express the ideal $I_{KKT}$ as an
intersection of super-ideals $I_{j,KKT}$. The super-ideal
$I_{j,KKT}$ denotes the ideal constructed by fixing a subset of
$p$ active constraints $\tilde g_i(u,x)$ among the set of all $q$
constraints $g_i(u,x)$ -- see (\ref{Eqn:gradideal}). The
corresponding Lagrange multipliers are denoted with $\tilde
\mu_i$. This leads to
\begin{equation}
\begin{array}{r}
I_{j,KKT}= \langle \; \nabla_{u}{J(u,x)} + \sum_{i=1}^{p}{\tilde
\mu_i \nabla_{u} \tilde g_i(u,x)}, \\
\tilde g_i(u,x) \; \rangle \\
\end{array}
\label{activeKKT}
\end{equation}
with the feasibility inequalities
\begin{equation}
\begin{array}{rcl}
\tilde\mu_i  & \geq & 0 \\
g_i(x,u) &\leq& 0\;.
\end{array}
\label{eqn::feas}
\end{equation}
Therefore, the ideal $I_{KKT}$ can be expressed as an intersection
of $\theta := \sharp(\{g_i(u,x)\}_{i = 1}^q) = 2^{q}$
super-ideals, where $\theta$ is the cardinality of the power set
of all $q$ constraints. Namely,
\begin{equation}
I_{KKT} = \bigcap_{j=1}^{\theta} I_{j,KKT} \; .
\label{eqn::bigcap}
\end{equation}
Relations (\ref{activeKKT}) and (\ref{eqn::feas}) lead to a large
number of super-ideals which are much better numerically
conditioned than the original problem, even though they are not
necessarily radical. Since many of the sub-varieties
$\mathcal{V}(I_{KKT})$ are empty, a Gr\"obner basis computation
for each ideal $I_{j,KKT}$ identifies these infeasible cases in
advance and reduces the subsequent companion matrix computations
tremendously by discarding them.

The number of solutions over $\bar{K}$ in the non-empty
sub-varieties $\VV_{j,KKT} = \VV(I_{j,KKT})$ can be calculated by
means of the Hilbert polynomial~(\cite{CoLO92}, Chapter 9, \S 3).
For zero-dimensional varieties this polynomial reduces to an
integer, which is equal to the number of solutions counting
multiplicity.

If the sub-variety has only a single solution, the coordinates
$u_i$ of the candidate solution can be computed analytically as a
\emph{rational} function of the parameters $x$. In this case, the
polynomials in the Gr\"obner basis from a set of linear equations
in the decision variables that can be solved analytically. For all
sub-varieties with more than one solution, a companion matrix has
to be computed. The result are companion matrices whose entries
are rational functions of the parameter $x$.

Specialization of the parameters gives a map from the field K to
the field R of real numbers. If the real parameters are chosen
generically enough, then the given Gr\"obner basis remains a
Gr\"obner basis, but for special choices of the parameters some
trouble may arise. For instance, it may happen that a
specialization leads to zero denominators. To handle this case,
comprehensive Gr\"obner bases can be used \cite{Weis92}. The
parametric computation is guaranteed to be correct only if the
sequence of leading coefficients of the result and the sequence of
greatest common denominators removed in the computations are
nonzero \cite{Weis92}. If ordinary methods such as Buchberger's
algorithm are used to compute Gr\"obner bases, these issues have
to be kept in mind.

A summary of the off-line algorithm appears in
Algorithm~\ref{offlinealgorithm}.
\begin{algorithm}
\caption{\emph{Off-line Part:}}
\begin{algorithmic}[1]
\REQUIRE Objective function $J(x,u)$ and constraints $g_i(x,u)
\leq 0$. \ENSURE Set of feasible sub-varieties $\VV_{j,KKT}$ with
their generalized companion matrices $M_{u_{j,i}}$ and
$M_{j,\tilde \mu_i}$, or an explicit function $u^*_{j,i}$ for
their candidate optimizer. \FORALL {combination of active and
inactive constraints} \STATE construct $I_{j,KKT}$ \STATE  calc.
Gr\"obner basis $G_j$ for $I_{j,KKT}$ \IF {$G_j = <1>$} \STATE
discard the super-ideal \ELSE \STATE calculate number of solutions
of $\VV_{j,KKT}$ by means of the Hilbert polynomial \IF {$\sharp
\VV_{j,KKT}=1$} \STATE Express all $u^*_{j,i}$ as rational
functions in the parameter $x$\ELSE
\STATE Compute generalized companion matrices $M_{j,u_{i}}$ and
$M_{j,\tilde \mu_i}$ for all decision variables $u_i$ \ENDIF
\ENDIF \ENDFOR \STATE \STATE \RETURN $M_{j,u_{i}}$ and
$M_{j,\tilde \mu_i}$, resp. $u^*_{j,i}$ and $\tilde \mu_{j,i}$
\end{algorithmic}
\label{offlinealgorithm}
\end{algorithm}

\subsection{On-line Part}
In order to evaluate the point coordinates of the KKT
sub-varieties, we need to compute eigenvectors and eigenvalues for
the companion matrices. Generally, eigenvalue computation cannot
be done parametrically. The parameter $x$ has to be fixed to a
numerical value and this computation is done on-line.

Given the precomputed generalized companion matrices $M_{j,u_i}$
and $M_{j,\tilde \mu_i}$ (resp. an explicit expression for all
sub-varieties with linear Gr\"obner basis) for all possible
feasible combinations of active and inactive constraints, the
on-line algorithm takes the value of the parameters $x$ to compute
the optimum $J^*$ and the optimizer $u^*$. The three main steps of
the algorithm are:
\begin{enumerate}
\item calculate all critical points \item remove infeasible
solutions \item find the feasible solution $u^*$ with the smallest
objective function value $J^*=J(u^*)$.
\end{enumerate}

Since all companion matrices have been computed parametrically,
the remaining part that has to be done is linear algebra. For
every non-empty sub-variety $\VV_{j,KKT}$, a set of right
eigenvectors $\{v\}$ is computed for the companion matrices
$M_{j,u_i}$ of the $j$-th sub-variety, see Theorem
\ref{th::eigenvalue}. Because all companion matrices for a
sub-variety $\VV_{j,KKT}$ commute pairwise, they form a
commutative sub-algebra within the non-commutative algebra of $l
\times l$ matrices, where $l$ is the companion matrix dimension
(\ref{eqn::basis}), see also \cite{CoLO98}. Therefore, it suffices
to calculate the eigenvectors for a single arbitrary matrix in
this sub-algebra, because they all share the same eigenvectors. To
avoid computational problems, we choose a matrix $M_{j,rand}$ in
this sub-algebra as a random linear combination of the companion
matrices associated with the decision variables $M_{j,u_i}$, i.e.
\begin{equation}
\begin{array}{ll}
\!\!\! M_{j,rand} &= c_1 M_{j,u_1}+ \cdots + \; c_{m}
M_{j,u_{m}} +
\\
& + \; c_{m+1} M_{j, \tilde \mu_1} + \cdots\ + \; c_{m+p}
M_{j,\tilde \mu_p } \; , \! \!
\end{array}
\end{equation}
where $c_i \in \R$ are randomly chosen. This ensures, with a low
probability of failure, that the corresponding eigenvalues will
all have algebraic multiplicity of one~(\cite{CoLO98}, Chapter 2,
\S 4).

The sets of eigenvectors $\{v\}_j$ can now be used to compute all
candidate critical points and their Lagrange multipliers $\tilde
\mu_{j,k}$ for the sub-variety $\VV_{j,KKT}$. To avoid unnecessary
computations, we first calculate the candidate Lagrange
multipliers $\tilde \mu_{j,i}$ for each sub-variety $\VV_{j,KKT}$.
In this way, complex or infeasible candidate points with
$\mu_{j,i} <0$ for some $i$ can be immediately discarded before the candidate
optimizers $u^*_{j,i}$ are computed. For all sub-varieties with
cardinality one, the problem of computing the critical points
reduces to an evaluation of the precomputed functions.

For all non-discarded candidate solutions, it remains to be
checked whether they are feasible, i.e. $g(u^*_{j,i},x_i) \leq 0$.
To achieve that, a set of feasible local candidate optimizers
$\mathcal{S} = \lbrace u^*_{j,i} \rbrace$ is initially calculated
by collecting all feasible candidate optimizers. After computing
the objective function value $J(u^*_{j,i},x)$ for all candidate
optimizers, the optimal solution
\begin{displaymath}
J^*=\min_{u^*_{j,i}\in \mathcal{S}} \limits J(u^*_{j,i},x)
\end{displaymath}
and the optimizer
\begin{equation*}
u^*_i=\argmin_{u^*_{j,i}\in \mathcal{S}} J(u^*_{j,i},x)
\end{equation*}
for the optimization problem~(\ref{Eqn:paramopt}) can be easily
obtained via discrete optimization over the finite set
$\mathcal{S}$.

A summary of the on-line algorithm can be seen in
algorithm~\ref{alg1}.

\begin{algorithm}
\caption{\emph{On-line Part:} Companion matrices $M_{u_i}$ and
$M_{ \tilde \mu_i}$ for all non-empty sub-varieties $\VV_{j,KKT}$,
resp. explicit expression for cardinality one sub-varieties has to
be provided.} \label{alg1}
\begin{algorithmic}[1]
\REQUIRE Value of the parameter $x$ (state measurement taken in
real time). \ENSURE Optimal cost $J^*$ and optimizer $u_{i}^*$.
\FORALL {feasible sub-varieties $\VV_{j,KKT}$ with $\sharp
\VV_{j,KKT}>1$} \STATE specialize parameter $x$ in $M_{u_i}$ and
$M_{\tilde \mu_i}$ \STATE  calc. a set of common eigenvectors
$\lbrace v \rbrace$ for the companion matrix $M_{j,rand}$ \STATE
solve $M_{j, \tilde \mu_i}  v = \tilde \mu_{j,i} v $ to obtain the
joint-eigenvalues, i.e. candidates for $\tilde \mu_{j,i}$ \STATE
discard all eigenvectors with corresp. $\tilde \mu_{j,i}<0$ \STATE
use the remaining eigenvectors to calc. joint-eigenvalues of
$M_{j,u_i}$ to obtain candidates for $u_{j,i}^*$ \ENDFOR
\FORALL{feasible sub-varieties $\VV_{j,KKT}$ with $\sharp
\VV_{j,KKT}=1$} \STATE evaluate $\tilde \mu_{j,i}(x)$ for all $i$
\IF {$\exists i : \;\tilde \mu_{j,i}(x)<0$} \STATE discard
sub-variety $\VV_{j,KKT}$\ELSE
\STATE evaluate $u^*_{j,i}(x)$\\
\ENDIF
\ENDFOR
\FORALL{evaluated candidate points $\lbrace u^*_{j,i} \rbrace_j$}
\IF {$g_k(u^*_{j,i},x) >0$}
\STATE discard candidate point $u^*_{j,i}$
\ELSE
\STATE evaluate $J(u^*_{j,i},x)$ \ENDIF
\ENDFOR

\STATE compare $J(u^*_{j,i},x)$ for the calculated candidates
$u^*_{j,i}$ and choose optimal $J^*$ and corresponding $u^*_i$
\STATE
\RETURN optimal cost $J^*$ and optimizer $u^*_{i}$
\end{algorithmic}
\label{onlinealgorithm}
\end{algorithm}

\section{OPTIMAL CONTROL APPLICATION}

In this section we fist give a description of the model predictive
control optimization problem to show the connection of parametric
optimization and optimal control.

\subsection{Nonlinear model predictive control}
Consider the nonlinear discrete-time system with state vector $x
\in \R^n$ and input vector $u \in \R^m$
\begin{equation}
{x(k+1) = f(x(k),u(k))} \label{eq::sys}
\end{equation}
subject to the inequality constraints
\begin{equation}
g(u(k),x(k)) \leq 0, \quad k=0,\ldots,N \; ,  \label{con::nl}
\end{equation}
where $N$ is the prediction horizon and $g\in$
$\R[x_1,\ldots,x_n,u_1,\ldots,u_m]^q$ is a vector polynomial
function representing the constraints of the problem. We consider
the problem of regulating system (\ref{eq::sys}) to the origin.
For that purpose, we define the following cost function
\begin{equation*}
\!\!\!J(U_0^{N-1},x_0)\! = \! \sum \limits
_{k=0}^{N-1}{L_{k}(x(k),u(k))} + L_{N}(x(N),u(N)) \; ,
\label{eq::cost}
\end{equation*}
where $U_0^{N-1} := [u(0), \ldots, u(N-1)]$ is the optimization
vector consisting of all the control inputs for $k=0,\ldots,N-1$
and $x(0)=x_0$ is the initial state of the system. Therefore,
computing the control input is equivalent to solving the following
nonlinear constrained optimization program
\begin{equation}
    \begin{array}{cl}
      &  \min_{u} \limits \, J(U_0^{N-1}, x_0)\\
&\quad \quad \quad x(k+1) = f(x(k),u(k))\\
&\quad \mbox{s.t.} \quad g(u(k),x(k))
\leq 0, \quad k=0,\ldots,N .
    \end{array}
    \label{eq::min}
\end{equation}

Forming a vector $u$ of decision variables with $u_k=u(k)$
and renaming $x(0)$, problem (\ref{eq::min}) is written in the
more compact form
\begin{equation}
\begin{array}{c}
        \min_{u} \limits \, J(u, x) \qquad \mbox{s.t.}
        \quad g(u,x) \leq 0,
    \end{array}
    \label{minprob}
\end{equation}
where $J(u,x)$ is a polynomial function in $u$ and $x$, $u \in
\R^m$ is the decision variable vector and the initial state
$x=x(0) \in \R^n$ is the parameter vector. This is exactly problem
(\ref{Eqn:paramopt}), a nonlinear parametric optimization problem.
Our goal is to obtain the vector of control moves $u$.

\subsection{Illustrative example}
In this section we illustrate the application of the proposed
method by means of a simple example. The off-line algorithm
including the algebraic methods and the case enumeration
(\ref{eqn::bigcap}) have been implemented in Maple. A
Maple-generated input file is used to initialize Matlab, in order
to compute the optimizer on-line.

Consider the Duffing oscillator \cite{JoSm87}, a nonlinear
oscillator of second order. An equation describing it in
continuous time is
\begin{equation}
\ddot y(t) + 2 \zeta \dot y(t) + y(t) + y(t)^3 = u(t),
\end{equation}
where $y \in \R$ is the continuous state variable and $u \in \R$
the control input. The parameter $\zeta$ is the damping
coefficient and is known (here $\zeta = 0.3$). The
control objective is to regulate the state to the origin. To
derive the discrete time model, forward difference approximation
is used (with a sampling period of $h = 0.05$ time units). The
resulting state space model with a discrete state vector $x \in
\R^2$ and input $u \in \R$ is
\begin{equation*}
\begin{array}{lll}
    \left[
    \begin{array}{c}
    x_1(k+1) \\
    x_2(k+1)
    \end{array}
    \right] &=&
    \left[
    \begin{array}{cc}
    1 & h \\
    -h & (1 - 2 \zeta h)
    \end{array}
    \right]
    \left[
    \begin{array}{c}
    x_1(k) \\
    x_2(k)
    \end{array}
    \right] \\[2ex]
    &+&
    \left[
    \begin{array}{c}
    0 \\
    h
    \end{array}
    \right] u(k)
    +
    \left[
    \begin{array}{c}
    0 \\
    -h x_1^3(k)
    \end{array}
    \right] .
\end{array}
\end{equation*}
An optimal control problem with prediction horizon $N=3$, weight
matrices
\begin{align*}
Q&=\left[\begin{matrix}1&0\\0&1
\end{matrix}\right],\\
R &= \frac{1}{10}
\end{align*}
and state-constraints
\begin{align*}
\Vert x(k+j) \Vert_{\infty} \leq 5 \;\;\forall j=1\dots N
\end{align*}
leads to the following optimization problem:
\begin{equation*}
\begin{array}{lll}
J^* &= &\\
&\min_{u(k),u(k+1),u(k+2)} \limits&\!\!\sum_{i=1}^3
\left[\begin{smallmatrix}x_1(k+i) &\!\!x_2(k+i)\end{smallmatrix}\right] Q
\left[\begin{smallmatrix}x_1(k+i)\\x_2(k+i)\end{smallmatrix}\right]\\
& &+
\sum_{i=0}^{2} u(k+i) R u(k+i)\\
\\
&\mbox{s.t.}\;\; \Vert x(k+j) \Vert_{\infty}&\!\! \leq 5
\;\;\forall j=1\dots N\,.
\end{array}
\end{equation*}
Of these twelve constraints there are ten constraints involving
$u(k+i)$, which have to be considered during the optimization. As
described in section~\ref{sec::algorithm} the KKT-variety will be
split in $2^{\sharp \{g_i\}} = 2^{10}=1024$ sub-varieties. For all
of them a Gr\"obner basis needs to be computed. It turns out that
only 29 of these are feasible, i.e. having a Gr\"obner basis
different from unity. Only these cases have to be further
considered in the online algorithm. Among them there are 24
sub-varieties $\mathcal{V}_{j,KKT}$ with a linear Gr\"obner basis.
For these, a closed form expression for the candidate optimizers
$u^*_{j,i}$ can be computed. For the remaining five cases
companion matrices have to be computed, requiring eigenvalue
computation in the on-line algorithm. These sub-varieties
$\mathcal{V}_{j,KKT}$ have five solutions counting multiplicities,
i.e. the companion matrices are $5 \times 5$ matrices.

\begin{figure}[ht]
\begin{center}
\input{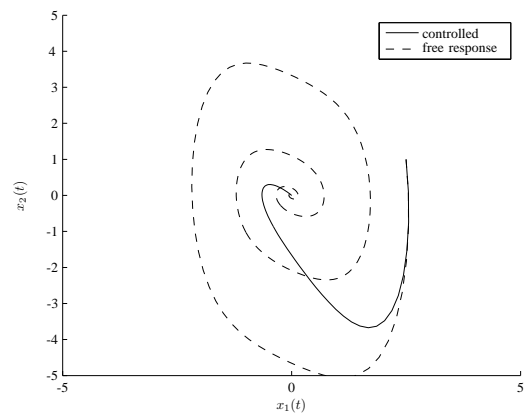}
\end{center}
\caption{State-space diagram of the Duffing oscillator}
\label{Fig:initial-response}
\end{figure}

\begin{figure}[ht]
\begin{center}
\input{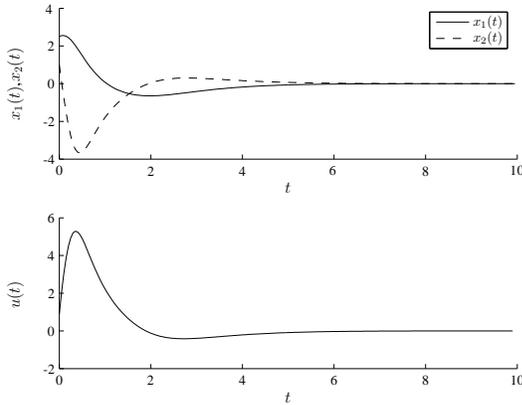}
\end{center}
\caption{State and input evolution of the controlled Duffing
oscillator} \label{Fig:State-Space}
\end{figure}

The trajectory of the controlled system starting from an initial
state of $x_1(0)=2.5 $ and $x_2(0)=1$ is shown in
Figure~\ref{Fig:initial-response}. Figure~\ref{Fig:State-Space}
shows the state-space evolution of the controlled Duffing
oscillator and its free response without the controller. In the
uncontrolled case, a weak dynamic behavior and a violation of the
constraint $x_2(t)>-5$ can be observed.

The precomputation of companion matrices and the solutions
$u^*_{j,i}$ took less than one minute on a Intel Pentium 3~GHz
with 1~GB RAM. The online algorithm needed less than 3.5~s to
obtain the global optimum even with a naive brute-force on-line
search algorithm for the minimization over the finite set of
candidate points. It has to be noted that most of the time of
these 3.5~s is consumed by the evaluation of expressions with the
Matlab Symbolic Math Toolbox. An efficient implementation, in C
for instance, would be orders of magnitude faster.

\section{CONCLUSIONS AND OUTLOOK}
The main contribution of this paper is a new algorithm for
nonlinear parametric optimization of polynomial functions subject
to polynomial constraints. The algorithm uses Gr\"obner bases and
the eigenvalue method for solving systems of polynomial equations,
to evaluate the map from the space of parameters to the
corresponding optimal value and optimizer. The algorithm is very
general, computationally robust and can be applied to a wide range
of problems.

The punchline of the proposed approach is the precomputation of
the generalized companion matrices, thus partially presolving the
optimization problem and moving the computational burden off-line.
The method has been developed with model predictive control in
mind. The connection to optimal control problems has been
illustrated by applying the method to the Duffing oscillator.

Finally, there is ongoing research on exploiting the structure of
specific control problems, including sparseness and genericity
assumption relaxation. More specifically, sparse resultant
techniques are investigated to compute the companion matrices.
Combining this method with recently proposed "Sum of Squares
Programming" methods, based on semi-definite representations of
finite varieties~\cite{Laur04}, seems to be a promising direction
for further research. Moreover, the integration of the proposed
scheme with dynamic programming is also explored.


\bibliographystyle{IEEEtran}
\bibliography{FoRS_ACC2006_bib}

\end{document}